\documentclass[12pt]{article}
\usepackage[latin1]{inputenc}
\usepackage{amssymb,amsmath,amsfonts}
\usepackage{xcolor}
\newtheorem{theorem}{Theorem}
\newtheorem{proposition}[theorem]{Proposition}
\newtheorem{lemma}[theorem]{Lemma}

\newtheorem{corollary}[theorem]{Corollary}

\newtheorem{remark}[theorem]{Remark}

\newtheorem{example}[theorem]{Example}

\newcommand{\R}{\mathbb{R}}
\newcommand{\Sp}{\mathbb{S}}
\newcommand{\U}{\mathcal{U}}
\renewcommand{\L}{\mathbb{L}}

\newcommand{\K}{{\cal K}}
\newcommand{\Ric}{{\cal R}}
\newcommand{\spa}{\mbox{span\,}}
\newcommand{\hess}{\mbox{Hess\,}}
\newcommand{\Ima}{\mbox{Im}}
\newcommand{\Hy}{\mathbb{H}}
\newcommand{\rank}{\mbox{rank }}

\newcommand{\grad}{\mbox{grad\,}}

\newcommand{\po}{{\hspace*{-1ex}}{\bf .  }}

\def\lp{{\langle\!\langle}}\vspace{2ex}
\def\rp{{\rangle\!\rangle}}
\def\<{{\langle}}
\def\>{{\rangle}}
\def\Sal{{\cal S}}

\def\a{\alpha}
\def\be{\begin{equation} }
\def\ee{\end{equation} }

\def\proof{\noindent\emph{Proof: }}
\def\qed{\ifhmode\unskip\nobreak\fi\ifmmode\ifinner
\else\hskip5 pt \fi\fi\hbox{\hskip5 pt \vrule width4 pt
height6 pt  depth1.5 pt \hskip 1pt }}
\newcommand\blfootnote[1]{
\begingroup
\renewcommand\thefootnote{}\footnote{#1}
\addtocounter{footnote}{-1}
\endgroup
}
\begin{document}

\title{Kaehler submanifolds of the real hyperbolic space}
\author{S. Chion and M. Dajczer}
\date{}
\maketitle

\begin{abstract}
The local classification of Kaehler submanifolds $M^{2n}$ of the 
hyperbolic space $\Hy^{2n+p}$ with low codimension $2\leq p\leq n-1$ 
under only intrinsic assumptions remains a wide open problem. The 
situation is quite different for submanifolds in the round sphere 
$\Sp^{2n+p}$, $2\leq p\leq n-1$, since  Florit, Hui and Zheng 
(\cite{FHZ}) have shown that the codimension has to be $p=n-1$ and 
then that any submanifold is just part of an extrinsic product of 
two-dimensional umbilical spheres in $\Sp^{3n-1}\subset\R^{3n}$. 
The main result of this paper is a version for Kaehler manifolds 
isometrically immersed into the hyperbolic ambient space of the 
result in \cite{FHZ} for spherical submanifolds. Besides, we 
generalize several results obtained by Dajczer and Vlachos 
(\cite{DV}). 
\end{abstract}
\blfootnote{\textup{2020} \textit{Mathematics Subject Classification}:
53B25,  53C40, 53C42.}

The study of isometric immersions of Kaehler manifolds 
$(M^{2n},J)$, $n\geq 2$, into spheres $\Sp^{2n+p}$ and hyperbolic 
spaces $\Hy^{2n+p}$  with low codimension $p$ was initiated by 
Ryan \cite{Ry}. He showed that for hypersurfaces there is only 
$M^4=\Sp^2\times\Sp^2\subset\Sp^5\subset\R^6$ in the sphere, 
and that in the hyperbolic space there is  
$M^4=\Hy^2\times\Sp^2\subset\Hy^5\subset\L^6$ besides the even 
dimensional horospheres. 
Later on Dajczer-Rodr\'iguez \cite{DR} proved that, regardless of 
the codimension, such isometric immersions do not occur if we 
require minimality.

The possibilities for Kaehler submanifolds with low codimension 
in spheres are rather restricted. In fact, it was shown by 
Florit-Hui-Zheng \cite{FHZ} that if we have an isometric immersion 
into the unit sphere 
$f\colon M^{2n}\to\Sp^{2n+p}_1$ with $p\leq n-1$ then $p=n-1$ 
and $f(M)\subset\Sp^{3n-1}_1\subset\R^{3n}=\R^3\times\cdots\times\R^3$ 
is an open subset of a Riemannian product of umbilical spheres 
$\{\Sp^2_{c_j}\}_{1\leq j\leq n}$ in $\R^3$ such that 
$1/c_1+\cdots+1/c_n=1$.

In \cite{FHZ} it was observed that for submanifolds in hyperbolic 
space a similar result as theirs is not possible due to the presence 
of the horospheres. In fact, there are the compositions $f=j\circ g$  
where $j\colon\R^{2n+p-1}\to\Hy^{2n+p}$ is a horosphere 
and $g\colon M^{2n}\to\R^{2n+p-1}$ an isometric immersion. 
In that respect, we observe that for any codimension there is 
an abundance of nonholomorphic Kaehler submanifolds, a class 
intensively studied in the last $25$ years with emphasis on 
the ones that are minimal. For an account of the basic 
facts on the subject of real Kaehler submanifolds we refer to 
Chapter $15$ in \cite{DT} and the references listed in 
\cite{DV}.
\vspace{1ex}

The main goal of this paper is to locally characterize the 
submanifolds described next.

\begin{example}\po\label{example} {\em Let the Kaehler manifold 
$M^{2n}$ be the Riemannian product of an hyperbolic plane and a 
set of two-dimensional round spheres given by
$$
M^{2n}=\Hy^2_{c_1}\times\Sp^2_{c_2}\times\cdots\times\Sp^2_{c_n}
\;\;\mbox{with}\;\; 1/c_1+\cdots +1/c_n=-1.
$$
Then let $f\colon M^{2n}\to\Hy^{3n-1}_{-1}$ be the submanifold
defined by $g=i\circ f\colon M^{2n}\to\L^{3n}$ where 
$i\colon\Hy^{3n-1}\to\L^{3n}$ is the inclusion into the flat 
Lorentzian space $\L^{3n}=\L^3\times\R^3\times\cdots\times\R^3$
and $g=g_1\times g_2\times\cdots\times g_n$ is the extrinsic 
product of umbilical surfaces  $g_1\colon \Hy^2_{c_1}\to\L^3$ 
and $g_j\colon\Sp^2_{c_j}\to\R^3$, $2\leq j\leq n$.
}\end{example}

To make our goal feasible it is necessary to remove the possibility 
for the submanifold to lay inside a horosphere. Such a task is fulfilled here inspired by the following sharp estimate given as 
Corollary $15.6$ in \cite{DT}. 
\vspace{2ex}

\noindent \emph{Let $f\colon M^{2n}\to\R^{2n+p}$, $1\leq p\leq n-1$, 
be an isometric immersion of a Kaehler manifold. At any $x\in M^{2n}$ 
there is a complex vector subspace $L^{2\ell}\subset T_xM$ with 
$\ell\geq n-p$ such that the sectional curvature of $M^{2n}$ satisfies 
$K(Z,JZ)\leq 0$ for any $Z\in L^{2\ell}$}. 
\vspace{1ex}

The following is the main result of this paper.

\begin{theorem}\label{classification}\po
Let $f\colon M^{2n}\to\Hy^{2n+p}_{-1}$, $p\leq n-1$, be 
an isometric immersion of a Kaehler manifold. Assume that at 
some point $x_0\in M^{2n}$ there is a complex vector subspace 
$V^{2m}\subset T_{x_0}M$ with $m\geq p$ such that the 
sectional curvature of $M^{2n}$ satisfies $K(S,JS)>0$ for 
any $0\neq S\in V^{2m}$. Then $p=n-1$ and $f(M)$ is an open subset 
of the submanifold given by Example \ref{example}.
\end{theorem}

The remaining of this paper is devoted to generalize several results 
due to Dajczer-Vlachos \cite{DV}. The following result improves their 
Theorem~$7$ that deals with the size of the dimension of the complex 
subspaces where the holomorphic sectional curvature is nonpositive. 
In addition, we establish an estimate for the Ricci curvature.
Moreover, in both cases the estimates now obtained are sharp.

\begin{theorem}\po\label{theorem}
Let $f\colon M^{2n}\to\Hy^{2n+p}$, $1\leq p\leq n-2$, be an isometric
immersion of a Kaehler manifold.  At any point $x\in M^{2n}$ there is 
a complex vector subspace $V^{2\ell}\subset T_xM$ with $\ell\geq n-p+1$ 
such that for any $S\in V^{2\ell}$ the sectional curvature of $M^{2n}$ 
satisfies $K(S,JS)\leq 0$ and the Ricci curvature that $Ric\,(S)\leq 0$.
\end{theorem}

For $p=1$ the above estimates follow trivially from the aforementioned 
result due to Ryan. For codimension $p=2$ they are a consequence of Theorem $1$ 
in \cite{DV} and, as is the case for $p=1$, with the stronger assertions 
$K(S,JS)=0=Ric\,(S)$. It is shown in \cite{DV} that locally $f=j\circ g$ 
is a composition as given above. To reach that conclusion one has to use 
that nonflat Kaehler hypersurfaces in Euclidean space, whose classification 
can be seen in \cite{DT} as Theorem~$15.14$, have only two nonzero simple 
principal curvatures. In particular, observe that when the hypersurface 
has a plane of positive sectional curvature these examples show that 
Theorem \ref{theorem} is sharp already for $p=2$. 
\vspace{1ex}

From the two results above we obtain the following generalization
of Theorem $3$ in \cite{DV} given there for codimension $p\leq n-2$.

\begin{corollary}\po A Kaehler manifold $M^{2n}$, $n\geq 2$, that at 
some point possesses positive holomorphic sectional curvature cannot 
be isometrically immersed in $\Hy^{3n-1}$. 
\end{corollary}

The next result was obtained in \cite{DV} for $p\leq n-2$ 
under the weaker hypothesis that the Omori-Yau weak maximum principle 
for the Hessian holds on $M^{2n}$. Under the assumptions that the 
Riemannian manifold is complete with sectional curvature bounded from 
below, we have from Theorem $2.3$ in \cite{AMR} that the Omori-Yau 
maximum principle for the Hessian holds.

\begin{theorem}\po\label{unbounded} Let $f\colon M^{2n}\to\Hy^{3n-1}$ 
be an isometric immersion of a complete Kaehler manifold with sectional 
curvature bounded from below. Then $f(M)$ is unbounded.
\end{theorem}

Finally, we consider the case of submanifolds with codimension two.

\begin{theorem}\po\label{m6} Let $f\colon M^{2n}\to\Hy^{2n+2}$, 
$n\geq 3$, be an isometric immersion of a Kaehler manifold which 
does not contain an open subset of flat points. Then either $n=3$ 
and $M^6\subset\Hy^2\times\Sp^2\times\Sp^2\subset\Hy^8\subset\L^9$
or there is a composition of isometric  immersion $f=j\circ g$ where 
$g\colon M^{2n}\to\R^{2n+1}$ is a real Kaehler hypersurface and 
$j\colon\R^{2n+1}\to\Hy^{2n+2}$ is a horosphere. 
\end{theorem}

In addition, the result extends Theorem $1$ in \cite{DV} by including the 
case $n=3$. Moreover, it generalizes this result as well as Theorem $6$ 
in \cite{DV} since now it is global and, with respect to the latter, it 
does not require to assume flat normal bundle.

\section{Several algebraic considerations}

\subsection{Some general facts}

Let $W^{p,p}$ denote a real vector space of dimension $2p$ 
endowed with an indefinite inner product of signature $(p,p)$. 
Thus $p$ is the maximal dimension of a vector 
subspace such that the induced inner product is either positive 
or negative definite. A vector subspace $L\subset W^{p,p}$ is 
called \emph{degenerate} if $L\cap L^\perp\neq 0$ and 
\emph{nondegenerate} if otherwise. Moreover, a vector subspace 
$L\neq 0$ is called \emph{isotropic} if it satisfies 
$L=L\cap L^\perp$.\vspace{1ex}

The following result is Sublemma $2.3$ in \cite{CD} and  
Corollary $4.3$ in \cite{DT}.

\begin{proposition}\po\label{l:decomposition}
Given a vector subspace $L\subset W^{p,p}$ there is a direct sum 
decomposition $W^{p,p}=\U\oplus\hat{\U}\oplus\mathcal{V}$ where 
$\U=L\cap L^\perp$ such that the vector subspace $\hat{\U}$ is 
isotropic, the vector subspace $\mathcal{V}=(\U\oplus\hat{\U})^\perp$ 
is nondegenerate and $L\subset\U\oplus\mathcal{V}$.
\end{proposition}

Let $V$ be a finite dimensional real vector space and 
$\varphi\colon V\times V\to W^{p,p}$ a bilinear 
form. Then $\varphi$ is called a \emph{flat bilinear form} if
$$
\<\varphi(X,Y),\varphi(Z,T)\>-\<\varphi(X,T),\varphi(Z,Y)\>=0
$$
for all $X,Z,Y,T\in V$. 
We denote the vector subspace generated by $\varphi$ by
$$
\mathcal{S}(\varphi)=\spa\{\varphi(X,Y)\colon X,Y\in V\}
$$
and say that $\varphi$ is \emph{surjective} if 
$\mathcal{S}(\varphi)=W^{p,p}$.
The (right) kernel $\varphi$ is defined by
$$
\mathcal{N}(\varphi)=\{Y\in V\colon\varphi(X,Y)=0
\;\mbox{for all}\;X\in V\}.
$$

If $V_1,V_2\subset V$ are  vector subspaces we denote 
$$
\Sal(\varphi|_{V_1\times V_2})
=\spa\{\varphi(X,Y)\colon \mbox{for all}\; X\in V_1\;\mbox{and}\;Y\in V_2\}.
$$

A vector $X\in V$ is called a (left) \emph{regular element} 
of $\varphi$ if $\dim\varphi_X(V)=r$ where
$$
r=\max\{\dim\varphi_X(V)\colon X\in V\}
$$
and $\varphi_X\colon V\to W^{p,p}$ is the linear 
transformation defined by 
$$
\varphi_XY=\varphi(X,Y).
$$
The set $RE(\varphi)$ of regular elements of 
$\varphi$ is easily seen to be an open dense subset 
of $V$; for instance see  Proposition $4.4$ in \cite{DT}. 

\begin{proposition}\po
Let $\varphi\colon V\times V\to W^{p,p}$ be a flat bilinear 
form. If $X\in RE(\varphi)$ then 
\be\label{moore}
\Sal(\varphi|_{V\times \ker\varphi_X})\subset
\varphi_X(V)\cap(\varphi_X(V))^\perp.
\ee 
\end{proposition}

\proof See Sublemma $2.4$ in \cite{CD} or Proposition $4.6$ in 
\cite{DT}.\qed

\subsection{A general bilinear form}
 
Let $V^{2n}$ and $\L^p$, $p\geq 2$, be real vector spaces 
such that there is $J\in Aut(V)$ which satisfies $J^2=-I$ 
and $\L^p$ is endowed with a Lorentzian inner product $\<\,,\,\>$. 
Then let $W^{p,p}=\L^p\oplus\L^p$ be endowed 
with the inner product of signature $(p,p)$ defined by
\be\label{metric}
\lp(\xi,\bar\xi),(\eta,\bar\eta)\rp
=\<\xi,\eta\>-\<\bar\xi,\bar\eta\>.
\ee

Let $\a\colon V^{2n}\times V^{2n}\to\L^p$ be 
a symmetric bilinear form and then let the bilinear 
form $\beta\colon V^{2n}\times V^{2n}\to W^{p,p}$ be defined by
\be\label{beta}
\beta(X,Y)=(\a(X,Y)+\a(JX,JY),\a(X,JY)-\a(JX,Y)).
\ee
Then $\beta$ satisfies 
\be\label{antihol}
\beta(X,JY)=-\beta(JX,Y)
\ee 
for any $X,Y\in V^{2n}$. Notice that $\beta(X,X)=(\zeta,0)$
for $X\in V^{2n}$ and $\zeta\in\L^p$.
Moreover, if  $\beta(X,Y)=(\xi,\eta)$ we have
\be\label{symmetries}
\beta(X,JY)=(\eta,-\xi),\;\;\beta(Y,X)=(\xi,-\eta)
\;\;\mbox{and}\;\;\beta(JY,X)=(\eta,\xi).
\ee 
In particular, it follows that $\mathcal{N}(\beta)$ is 
a $J$-invariant vector subspace of $V^{2n}$.
\vspace{1ex}

In the sequel $U_0^s\subset\L^p$ is the $s$-dimensional 
vector subspace $U_0^s=\pi_1(\Sal(\beta))$
where $\pi_1\colon W^{p,p}\to\L^p$ denotes the projection 
onto the first component of $W^{p,p}$. Hence
$$
U_0^s=\spa\{\a(X,Y)+\a(JX,JY)\colon X,Y\in V^{2n}\}.
$$
 
\begin{proposition}\label{betaimage}\po
The bilinear form $\beta\colon V^{2n}\times V^{2n}\to W^{p,p}$
satisfies
\be\label{ebetaimage}
\mathcal{S}(\beta)=U_0^s\oplus U_0^s.
\ee 
Moreover,  if the vector subspace $\Sal(\beta)$ is 
degenerate then $1\leq s\leq p-1$ and there is a nonzero
light-like vector $v\in U_0^s$ such that 
\be\label{uimage}
\Sal(\beta)\cap(\Sal(\beta))^\perp=\spa\{v\}\oplus\spa\{v\}.
\ee  
\end{proposition}
 
\proof If $\sum_{j=1}^k\beta(X_j,Y_j)=(\xi,\eta)$ 
we obtain from \eqref{symmetries} that
$$
\sum_{j=1}^k\beta(Y_j,X_j)=(\xi,-\eta),\;\;
\sum_{j=1}^k\beta(X_j,JY_j)=(\eta,-\xi)\;\;\mbox{and}\;\;
\sum_{j=1}^k\beta(JY_j,X_j)=(\eta,\xi).
$$
Hence $(\xi,0),(0,\xi),(\eta,0)\in\Sal(\beta)$ and 
thus $\mathcal{S}(\beta)\subset U_0^s\oplus U_0^s$.
If $(\xi,\eta)\in U_0^s\oplus U_0^s$  
there are $\delta,\bar{\delta}\in \L^p$ such that 
$(\xi,\delta),(\eta,\bar{\delta})\in\Sal(\beta)$, and 
as just seen $(\xi,\eta)\in\Sal(\beta)$. 

Since $\Sal(\beta)\neq 0$ and \eqref{ebetaimage}
is satisfied then if $\Sal(\beta)$ is a degenerate subspace we have 
that $1\leq s\leq p-1$.  If $\U=\Sal(\beta)\cap(\Sal(\beta))^\perp$ 
we claim that $\U=U_1\oplus U_1$ where $U_1=\pi_1(\U)$. 
We have from \eqref{symmetries} that
$$ 
\lp\beta(X,Y),(\eta,-\xi)\rp=\lp\beta(X,JY),(\xi,\eta)\rp,
\;\;\lp\beta(X,Y),(\xi,-\eta)\rp=\lp\beta(Y,X),(\xi,\eta)\rp
$$
and 
$$
\lp\beta(X,Y),(\eta,\xi)\rp=-\lp\beta(JY,X),(\xi,\eta)\rp.
$$
Hence if $(\xi,\eta)\in\U$ then also $(\eta,-\xi)$,
$(\xi,-\eta)$, $(\eta,\xi)\in\U$. Thus
$(\xi,0),(0,\xi),(\eta,0)\in\U$ and hence
$\U\subset U_1\oplus U_1$. If $(\xi,\eta)\in U_1\oplus U_1$ 
there are $\delta,\bar{\delta}\in \L^p$ such that 
$(\xi,\delta),(\eta,\bar{\delta})\in\U$ and thus $(\xi,\eta)\in\U$ 
proving the claim. It follows from \eqref{metric} and 
the claim that the vector subspace $U_1\subset\L^p$ 
is isotropic and thus \eqref{uimage} holds.
\vspace{2ex}\qed

Given $X\in V^{2n}$ we denote $N(X)=\ker B_X$ where $B_X=\beta_X$. 
It follows from \eqref{symmetries} that the vector subspace
$N(X)$ is $J$-invariant. 

\begin{lemma}\label{charac}\po Let the bilinear form
$\beta$ be flat and $X\in RE(\beta)$ satisfy 
$\beta|_{V\times N(X)}\neq 0$.
Then there is a nonzero light-like vector $v\in\L^p$ such that
\be\label{lightlike}
\spa\{v\}\oplus\spa\{v\}\subset\Sal(\beta|_{V\times N(X)})\subset
B_X(V)\cap (B_X(V))^\perp.
\ee
Moreover, if $(v',w')\in B_X(V)$ where $v',w'$ are 
light-like vectors then $v',w'\in\spa\{v\}$. 
\end{lemma}

\proof First notice that the second inclusion in \eqref{lightlike} 
is just \eqref{moore}. If 
$\beta|_{N(X)\times N(X)}\neq 0$ then the vector subspace
$\Sal(\beta|_{N(X)\times N(X)})$ is isotropic. Then by
\eqref{uimage} there is $v\in\L^p$ such that 
$\spa\{v\}\oplus\spa\{v\}=\Sal(\beta|_{N(X)\times N(X)})$. 
Suppose that $\beta|_{N(X)\times N(X)}=0$. By assumption 
there are $Y\in V^{2n}$ and $Z\in N(X)$ such that 
$\beta(Y,Z)=(v,w)\neq 0$. Since the vector subspace
$\Sal(\beta|_{V\times N(X)})$ is isotropic we have
$$
0=\lp\beta(Y,Z),\beta(Y,Z)\rp=\|v\|^2-\|w\|^2
$$
whereas from \eqref{symmetries} and the flatness of 
$\beta$ we obtain 
$$
0=\lp\beta(Y,Z),\beta(Z,Y)\rp=\|v\|^2+\|w\|^2.
$$
Thus the vectors $v,w\in\L^p$ are both light-like. 

It suffices to argue for $v\neq 0$ since $\beta(Y,JZ)=(w,-v)$. 
Since $N(X)$ is $J$-invariant  and $\Sal(\beta|_{V\times N(X)})$ 
is  isotropic we obtain using \eqref{symmetries} that
$$
0=\lp\beta(Y,Z),\beta(Y,JZ)\rp=2\<v,w\>
$$
and hence $w=av$. Then from
$$
\beta(Y,Z+aJZ)=(a^2+1)(v,0)
\;\;\mbox{and}\;\;\beta(Y,aZ-JZ)=(a^2+1)(0,v)
$$
we obtain the first inclusion in \eqref{lightlike}.

Let $B_XZ=(v',w')$ be as in the statement. By \eqref{lightlike} 
we have
$$
\<v,v'\>=\lp(v,0),(v',w')\rp=\lp(v,0),\beta(X,Z)\rp=0
$$
and thus $v'\in\spa\{v\}$. Since $B_XJZ=(w',-v')$ then also
$w'\in\spa\{v\}$.\qed 

\begin{proposition}\label{mainlemma}\po Let the bilinear form
$\beta\colon V^{2n}\times V^{2n}\to W^{p,p}$, $p\leq n$, be 
surjective and flat. Then we have
\be\label{estimnucleo2}
\dim\mathcal{N}(\beta)\geq 2n-2p.
\ee
Moreover, if $\mathcal{N}(\beta)=0$ and $X\in RE(\beta)$
then $B_X\colon V^{2n}\to W^{p,p}$ is an isomorphism.
\end{proposition}

\proof  If $X\in RE(\beta)$ we have $\dim N(X)\geq 2n-2p$.
Since $\mathcal{N}(\beta)\subset N(X)$, if we have that 
$N(X)=0$ for some $X\in RE(\beta)$ then the result
holds trivially. Thus, it remains to argue when $N(X)\neq 0$ 
for any $X\in RE(\beta)$. In this case, we also show that
$\mathcal{N}(\beta)\neq 0$ and this gives the second statement.

If $\beta|_{V\times N(X)}=0$ for some $X\in RE(\beta)$ then 
$N(X)=\mathcal{N}(\beta)$. Then $\mathcal{N}(\beta)\neq 0$
and $\dim\mathcal{N}(\beta)=\dim N(X)\geq 2n-2p$.  Hence, 
we assume that $\beta|_{V\times N(X)}\neq 0$
for any $X\in RE(\beta)$. Fix a vector $X\in RE(\beta)$. 
By Lemma \ref{charac} there is a nonzero light-like vector 
$v\in\L^p$ such that
$$
\spa\{v\}\oplus\spa\{v\}\subset \U^\tau(X)
=B_X(V)\cap (B_X(V))^\perp.
$$
On one hand, since $\Sal(\beta)=W^{p,p}$ and the subset 
$RE(\beta)$ is dense there are $Y\in RE(\beta)$ and $
Z\in V^{2n}$ such that $\lp\beta(Y,Z),(v,0)\rp\neq 0$. 
On the other hand, since $Y\in RE(\beta)$ then Lemma~\ref{charac} 
yields a nonzero light-like vector $w\in\L^p$ such that
$\spa\{w\}\oplus\spa\{w\}\subset\Sal(\beta|_{V\times N(Y)})
\subset B_Y(V)\cap(B_Y(V))^\perp$, and hence
$\lp\beta(Y,Z),(w,0)\rp=0$. Thus the vectors $v$
and $w$ are linearly independent.

From \eqref{moore} we obtain that $B_Y(N(X))\subset\U^\tau(X)$.
If $(av,bv)\in B_Y(N(X))$ then the second part of Lemma \ref{charac} 
gives that $av,bv\in\spa\{w\}$, and hence $a=b=0$. Thus
$\spa\{v\}\oplus\spa\{v\}\cap B_Y(N(X))=0$. Consequently, 
for $B_Y|_{N(X)}\colon N(X)\to\U^\tau(X)$  and since 
$\spa\{v\}\oplus\spa\{v\}\subset \U^\tau(X)$, we have that 
$N_1=\ker B_Y|_{N(X)}$ satisfies
\be\label{n1}
\dim N_1\geq\dim N(X)-\tau+2.
\ee

Proposition \ref{l:decomposition} applied to
$B_X(V)\subset W^{p,p}$ yields a decomposition
$$
W^{p,p}
=\U^\tau(X)\oplus\hat\U^\tau(X)\oplus\mathcal{V}^{p-\tau,p-\tau}
$$
verifying that $B_X(V)\subset\U(X)\oplus\mathcal{V}$ among 
other properties. Thus $\dim B_X(V)\leq 2p-\tau$ and hence 
$\dim N(X)\geq 2n-2p+\tau$. It follows from \eqref{n1} that 
$\dim N_1\geq 2n-2p+2\geq 2$.

We prove that $N_1=\mathcal{N}(\beta)$ which gives
$\dim\mathcal{N}(\beta)\geq 2n-2p+2>0$ that is even a better
estimate 
than \eqref{estimnucleo2}. 
Since $N_1=N(X)\cap N(Y)$ by \eqref{moore} then
$\Sal(\beta|_{N_1\times N_1})$ is an isotropic vector subspace
unless $\beta|_{N_1\times N_1}=0$.  In the former case, we have
from Proposition \ref{betaimage} applied to $\beta|_{N_1\times N_1}$
that there is a nonzero light-like vector $z\in\L^p$ such that
$$
\Sal(\beta|_{N_1\times N_1})=
\Sal(\beta|_{N_1\times N_1})\cap(\Sal(\beta|_{N_1\times N_1}))^\perp
=\spa\{z\}\oplus\spa\{z\}.
$$
Since $N_1\subset N(X)$ we obtain from \eqref{lightlike} that
$(z,0)\in B_X(V)$ and, similarly, we have that $(z,0)\in B_Y(V)$. 
Then Lemma \ref{charac} yields $z\in\spa\{v\}\cap\spa\{w\}=0$,
which is not possible. We conclude that $\beta|_{N_1\times N_1}=0$.

If $\beta|_{V\times N_1}\neq 0$ there are vectors $Z\in V^{2n}$ and
$T\in N_1$ such that $\beta(Z,T)=(\xi,\eta)\neq 0$.
Then \eqref{symmetries} and the flatness of $\beta$ give
$$
0=\lp\beta(Z,T),\beta(T,Z)\rp=\|\xi\|^2+\|\eta\|^2.
$$
By \eqref{moore} the vector subspace
$\Sal(\beta|_{V\times N_1})$ is isotropic and thus
$$
0=\lp\beta(Z,T),\beta(Z,T)\rp=\|\xi\|^2-\|\eta\|^2.
$$
Thus the vectors  $\xi,\eta$ are light-like. Then 
from \eqref{moore} and the second part of Lemma \ref{charac}
we obtain that $\xi,\eta\in\spa\{v\}\cap\spa\{w\}=0$, and this 
is a contradiction. Then $\beta|_{V\times N_1}=0$ as wished.\qed

\subsection{A special flat bilinear form}

Throughout this section $V^{2n}$ is endowed with a 
positive definite inner product denoted by $(\,,\,)$ with respect 
to which $J\in Aut(V)$ is an isometry. In addition, we assume that 
there exists a time-like of unit length vector $w\in\L^p$ such that
\be\label{shapeid}
\<\a(X,Y),w\>=-(X,Y)
\ee
for any $X,Y\in V^{2n}$.
\vspace{1ex}

Under the assumptions above $\beta\colon V^{2n}\times V^{2n}\to W^{p,p}$ 
given by \eqref{beta} satisfies  
\be\label{conditionalpha}
\lp\beta(X,Y),(w,0)\rp=\<\a(X,Y)+\a(JX,JY),w\>=-2(X,Y).
\ee 
In particular, if $\beta(X,Y)=0$ then $(X,Y)=0$ and  
\be\label{nbeta}
\beta(X,X)\neq 0\;\;\mbox{if}\;\; 0\neq X\in V^{2n}\;\;\mbox{thus}\;\;
\mathcal{N}(\beta)\neq 0.
\ee

\begin{proposition}\po\label{decompor}
Let the bilinear form $\beta\colon V^{2n}\times V^{2n}\to W^{p,p}$ 
be flat and the vector subspace $\Sal(\beta)$ degenerate.
If $v\in U_0^s$ is a light-like vector as in \eqref{uimage} 
then the plane $L=\spa\{v,w\}\subset\L^p$ is Lorentzian. Moreover, 
choosing $v$  such that $\<v,w\>=-1$ and setting
$\beta_1=\pi_{L^\perp\times L^\perp}\circ\beta$, we have
\be\label{betadecomp}
\beta(X,Y)=\beta_1(X,Y)+2((X,Y)v,(X,JY)v)
\ee 
for any $X,Y\in V^{2n}$. Furthermore, if $s\leq n$ then
\be\label{estimnucleo}
\dim\mathcal{N}(\beta_1)\geq 2n-2s+2.
\ee
\end{proposition}

\proof The plane $L$ is trivially  Lorentzian. We choose $v$ 
such that $v=u+w$ where $u$ is a space-like unit vector orthogonal 
to $w$. Then $L=\spa\{u ,w\}$.

We have from \eqref{uimage} that
\be\label{vdeg}
0=\lp\beta(X,Y),(v,0)\rp=\<\a(X,Y)+\a(JX,JY),v\>
\ee 
for any $X,Y\in V^{2n}$. Then from \eqref{conditionalpha} and 
\eqref{vdeg} we have
$$
\<\a(X,Y)+\a(JX,JY),u\>=2(X,Y)
$$
for any $X,Y\in V^{2n}$.  From \eqref{conditionalpha} we obtain
\begin{align}\label{betacomp}
\a(X,Y)+\a(JX,JY)&=\a_{L^\perp}(X,Y)+\a_{L^\perp}(JX,JY)
+\<\a(X,Y)+\a(JX,JY),u\>u\nonumber\\
&\quad-\<\a(X,Y)+\a(JX,JY),w\>w\nonumber\\
&=\a_{L^\perp}(X,Y)+\a_{L^\perp}(JX,JY)+2(X,Y)v
\end{align}
and then \eqref{beta} gives \eqref{betadecomp}.

We have seen that $v$ can be chosen so that $\<v,w\>=-1$ and 
that \eqref{betadecomp} holds. From \eqref{ebetaimage} and 
\eqref{uimage} we obtain that $w\notin U_0^s$. Hence 
$\dim(U_0^s+L^\perp)=p-1$.  Then from
$$
\dim(U_0^s+L^\perp)=\dim U_0^s+\dim L^\perp-\dim U_0^s\cap L^\perp
$$
we have that $U_1=U_0^s\cap L^\perp$ satisfies 
\be\label{dim}
\dim U_1=s-1.
\ee
Hence $\Sal(\beta_1)=U_1^{s-1}\oplus U_1^{s-1}$ from \eqref{ebetaimage}, 
\eqref{uimage} and \eqref{betadecomp}.

From \eqref{betadecomp} the  bilinear form
$\beta_1\colon V^{2n}\times V^{2n}\to L^\perp\oplus L^\perp$
is flat. Let $X\in RE(\beta_1)$ and set 
$N_1(X)=\ker B_{1X}$ where $B_{1X}Y=\beta_1(X,Y)$. To
obtain \eqref{estimnucleo} it suffices to show that 
$N_1(X)=\mathcal{N}(\beta_1)$ since then 
$\dim\mathcal{N}(\beta_1)=\dim N_1(X)\geq 2n-2\dim U_1$.

Let $\beta_1(Y,Z)=(\xi,\eta)$ where $Y,Z\in N_1(X)$. 
From \eqref{moore} and \eqref{symmetries} we have
$$
0=\lp\beta_1(Y,Z),\beta_1(Z,Y)\rp
=\lp(\xi,\eta),(\xi,-\eta)\rp=\|\xi\|^2+\|\eta\|^2.
$$
Thus $\beta_1|_{N_1(X)\times N_1(X)}=0$ since the inner product 
induced on $U_1$ is positive definite. Now let 
$\beta_1(Y,Z)=(\delta,\zeta)$ where $Y\in V^{2n}$ and $Z\in N_1(X)$. 
The flatness of $\beta_1$ and \eqref{symmetries} yield
$$
0=\lp\beta_1(Y,Z),\beta_1(Z,Y)\rp
=\lp(\delta,\zeta),(\delta,-\zeta)\rp
=\|\delta\|^2+\|\zeta\|^2
$$
and hence $\beta_1|_{V\times N_1(X)}=0$.\qed

\begin{proposition}\po\label{umbilical}
Let the bilinear form $\beta\colon V^{2n}\times V^{2n}\to W^{p,p}$ 
be flat and the vector subspace $\Sal(\beta)$ degenerate.
If $\gamma\colon V^{2n}\times V^{2n}\to W^{p,p}$ 
is the bilinear form defined by 
\be\label{gamma}
\gamma(X,Y)=(\a(X,Y),\a(X,JY))
\ee
assume that
\be\label{product}
\lp\beta(X,Y),\gamma(Z,T)\rp=\lp\beta(X,T),\gamma(Z,Y)\rp
\ee 
for any $X,Y,Z,T\in V^{2n}$. If $s\leq n-1$ we have the following:
\begin{itemize}
\item[(i)] For $v\in U_0^s$ satisfying \eqref{uimage} 
we have
\be\label{umbilic}
\<\a(X,Y),v\>=0\;\;\mbox{for any}\;\; X,Y\in V^{2n}.
\ee 
\item[(ii)] Choosing $v\in U_0^s$ satisfying \eqref{uimage} such 
that $\<v,w\>=-1$ then with respect to  the Lorentzian
plane $L=\spa\{v,w\}$ we have 
\be\label{alphapar}
\a(X,Y)=\a_{L^\perp}(X,Y)+(X,Y)v\;\;\mbox{for any}\;\; X,Y\in V^{2n}.
\ee  
\item[(iii)] There is a $J$-invariant vector subspace 
$P^{2m}\subset V^{2n}$ with $m\geq n-s+1$ such that for any 
$S\in P^{2m}$ we have
$$
\K(S)=\<\a(S,S),\a(JS,JS)\>-\|\a(S,JS)\|^2\leq 0
$$ 
and 
$$
\Ric(S)=\sum_{i=1}^{2n}(\<\a(X_i,X_i),\a(S,S)\>-\|\a(X_i,S)\|^2)\leq 0
$$
where $\{X_i\}_{1\leq i\leq 2n}$ is an orthonormal basis of $V^{2n}$.
\end{itemize}
\end{proposition}

\proof $(i)$ It suffices to argue for $v\in U_0^s$ such that
$\<v,w\>=-1$.
From \eqref{betadecomp} and \eqref{gamma} we obtain
\be\label{product2}
\lp\gamma(X,Y),\beta(S,S)\rp=2\<\a(X,Y),v\>
\ee
for any $S\in\mathcal{N}(\beta_1)$ of unit length and 
$X,Y\in V^{2n}$.  Since we have from \eqref{symmetries} and 
\eqref{betadecomp} that $\beta(S,Y)=\beta(Y,S)=0$ for any 
$S\in\mathcal{N}(\beta_1)$ and $Y\in\{S,JS\}^\perp$, then  
\eqref{product} and \eqref{product2} give 
\be\label{follows}
\<\a(X,Y),v\>=0
\ee
for any $X\in V^{2n}$ and $Y\in\{S,JS\}^\perp$ where 
$S\in\mathcal{N}(\beta_1)$. Since $s\leq n-1$ then  
\eqref{estimnucleo} gives 
$\dim\mathcal{N}(\beta_1)\geq 4$ and now 
\eqref{umbilic} follows from \eqref{follows}.
\vspace{1ex}

\noindent $(ii)$ From Proposition \ref{decompor} the plane 
$L=\spa\{v,w\}$ is Lorentzian.  There is $u\perp w$ space-like 
of unit length such that $v=u+w$. Hence \eqref{shapeid} and 
\eqref{umbilic} give $\<\a(X,Y),u\>=(X,Y)$. Now since
$$
\a(X,Y)=\a_{L^\perp}(X,Y)+\<\a(X,Y),u\>u-\<\a(X,Y),w\>w
$$
then \eqref{alphapar} follows from \eqref{shapeid}.
\vspace{1ex}

\noindent $(iii)$ We choose $v\in U_0^s$ satisfying \eqref{uimage} 
such that $\<v,w\>=-1$. From \eqref{alphapar} we have
$$
\gamma(X,Y)=
(\a_{L^\perp}(X,Y)+(X,Y)v,\a_{L^\perp}(X,JY)+(X,JY)v).
$$
Set $P^{2m}=\mathcal{N}(\beta_1)$. From \eqref{estimnucleo} we have 
$2m=\dim\mathcal{N}(\beta_1)\geq 2n-2s+2$.
It follows from \eqref{betadecomp} that
$\beta(Z,S)=2((Z,S)v,(Z,JS)v)$ for any $S\in P^{2m}$
and $Z\in V^{2n}$. Then \eqref{product} gives
$$
\lp\gamma(X,S),\beta(Z,Y)\rp=\lp\gamma(X,Y),\beta(Z,S)\rp=0
$$
for any $S\in P^{2m}$ and $X,Y,Z\in V^{2n}$.
Thus the vector subspaces $\Sal(\gamma|_{V\times P})$ 
and $\Sal(\beta)$ are orthogonal.  From \eqref{dim}
we have 
\be\label{deco}
U_0^s=U_1^{s-1}\oplus\spa\{v\}
\ee
where $U_1^{s-1}=U_0^s\cap L^\perp$.
Then by \eqref{ebetaimage} the vector subspaces
$\Sal(\gamma|_{V\times P})$ and $U_1^{s-1}\oplus U_1^{s-1}$ are 
orthogonal and therefore
$$
\<\a(X,S),\xi\>=\lp\gamma(X,S),(\xi,0)\rp=0
$$
for any $X\in V^{2n}$, $S\in P^{2m}$ and $\xi\in U^{s-1}_1$.
Since $U_1^{s-1}\subset L^\perp$ then 
\be\label{nucleoa}
\a_{U_1}(X,S)=0
\ee 
for any $X\in V^{2n}$ and $S\in P^{2m}$.

Let $U_2^{p-s-1}\subset L^\perp$ be given by the orthogonal 
decomposition
$\L^p=U_1^{s-1}\oplus U_2^{p-s-1}\oplus L$.
By \eqref{ebetaimage} and \eqref{deco} we have 
$$
\<\a(X,Y)+\a(JX,JY),\xi_2\>
=\lp\beta(X,Y),(\xi_2,0)\rp=0
$$
for any $X,Y\in V^{2n}$ and $\xi_2\in U_2^{p-s-1}$.  Thus 
\be\label{pluri}
\a_{U_2}(X,Y)=-\a_{U_2}(JX,JY)
\ee 
for any $X,Y\in V^{2n}$.

From \eqref{alphapar}, \eqref{nucleoa}, \eqref{pluri} and 
since the inner product induced on $U_2^{p-s-1}$ is positive 
definite, we have
$$
\K(S)=\<\a_{U_2}(S,S),\a_{U_2}(JS,JS)\>-\|\a_{U_2}(S,JS)\|^2
=-\|\a_{U_2}(S,S)\|^2-\|\a_{U_2}(S,JS)\|^2\leq 0
$$
for any $S\in P^{2m}$. Also
\begin{align*}
\Ric(S)&=\sum_{i=1}^n(\<\a_{U_2}(Y_i,Y_i),\a_{U_2}(S,S)\>
-\|\a_{U_2}(Y_i,S)\|^2+\<\a_{U_2}(JY_i,JY_i),\a_{U_2}(S,S)\>\\
&\hspace{5ex}-\|\a_{U_2}(JY_i,S)\|^2)\\
&=-\sum_{i=1}^n(\|\a_{U_2}(Y_i,S)\|^2
+\|\a_{U_2}(JY_i,S)\|^2)\leq 0
\end{align*}
for any $S\in P^{2m}$ and an orthonormal basis 
$\{Y_j,JY_j\}_{1\leq j\leq n}$.\qed

\begin{lemma}\label{kercod2}\po Let the bilinear form
$\beta\colon V^{2n}\times V^{2n}\to W^{p,p}$ be flat. 
Assume that the vector subspace $\Sal(\beta)$ is nondegenerate.
If $Z\in V^{2n}$ satisfies $\dim N(Z)=2n-2$ then
$B_ZZ=(\xi,0)\neq 0$, the vector subspace
$B_Z(V)=\spa\{\xi\}\oplus\spa\{\xi\}$ is nondegenerate and
the decomposition
\be\label{kercod2decomp}
\Sal(\beta)=B_Z(V)\oplus\Sal(\beta|_{N(Z)\times N(Z)})
\ee
is orthogonal.
\end{lemma}

\proof We have $0\neq B_ZZ=(\xi,0)$ and $B_ZJZ=-(0,\xi)$.
If $Y\in N(Z)$ then
$$
0=\lp B_ZY,(w,0)\rp=-2(Z,Y)\;\;\mbox{and}\;\;
0=\lp B_ZY,(0,w)\rp=-2(JZ,Y)
$$
and thus $N(Z)\subset\{Z,JZ\}^\perp$.
Since $\dim N(Z)=2n-2$ then $N(Z)=\{Z,JZ\}^\perp$ and
$\dim B_Z(V)=2$ gives $B_Z(V)=\spa\{\xi\}\oplus\spa\{\xi\}$. 
From \eqref{antihol} and \eqref{symmetries} we obtain
$$
\Sal(\beta)
=B_Z(V)+\Sal(\beta|_{N(Z)\times N(Z)}).
$$
The flatness of $\beta$ gives
$$
\lp\beta(Z,Y),\beta(S,T)\rp=\lp B_ZT,\beta(S,Y)\rp=0
$$
for any $S,T\in N(Z)$ and $Y\in V^{2n}$. Hence the
decomposition \eqref {kercod2decomp} is orthogonal and
thus the vector subspace $B_Z(V)$ is nondegenerate.
\vspace{2ex}\qed

The following is the fundamental result obtained in this first 
section although it is not needed for the proof 
of Theorem \ref{theorem}. We observe that Lemma $7$ in \cite{FHZ} 
is a Riemannian version of the result. For an alternative 
version of the Riemannian result in the spirit of this paper see 
Proposition $2.6$ in \cite{CDa}. 

\begin{proposition}\label{diag}\po Let the bilinear form 
$\beta\colon V^{2n}\times V^{2n}\to W^{p,p}$, $s\leq n$, be flat.
Assume that the subspace $\Sal(\beta)$ is nondegenerate and  
\eqref{product} holds. For $p\geq 4$ assume further that there is 
no nontrivial $J$-invariant vector subspace $V_1\subset V^{2n}$ 
such that the subspace $\Sal(\beta|_{V_1\times V_1})$ is degenerate 
and that $\dim\Sal(\beta|_{V_1\times V_1})\leq\dim V_1-2$.
Then $s=n$ and there is an orthogonal basis $\{X_i,JX_i\}_{1\leq i\leq n}$ 
of $V^{2n}$ such that:
\begin{itemize}
\item[(i)] $\beta(Y_i,Y_j)=0\;\mbox{if}\; i\neq j
\;\mbox{and}\;Y_k\in\spa\{X_k,JX_k\}$\,\mbox{for}\,k=i,j. 
\item[(ii)] The vectors 
$\{\beta(X_j,X_j),\beta(X_j,JX_j)\}_{1\leq j\leq n}$ 
form an orthonormal basis of $\Sal(\beta)$.
\end{itemize}
\end{proposition}

\proof 
Since $\mathcal{N}(\beta)=0$ then by Proposition \ref{mainlemma} 
we have that $s=n$.  Moreover, we have that 
$B_X\colon V^{2n}\to U_0^n\oplus U_0^n\subset W^{p,p}$ is an isomorphism
for any $X\in RE(\beta)$.
\vspace{1ex}

\noindent{\it Fact} $1$. If $n\geq 2$ then there are nonzero vectors 
$X,Y\in V^{2n}$ satisfying $\beta(X,Y)=0$. 
\vspace{1ex}

We have that $RE(\beta)$ is open and dense in $V^{2n}$ and 
that $B_X\colon V^{2n}\to U_0^n\oplus U_0^n$ is an isomorphism for any 
$X\in RE(\beta)$. Hence, there is a basis $Z_1,\ldots,Z_{2n}$ of 
$V^{2n}$ with $Z_2\not\in\spa\{Z_1,JZ_1\}$ such that  
$\{\beta(Z_1,Z_j)\}_{1\leq j\leq 2n}$ and
$\{\beta(Z_2,Z_j)\}_{1\leq j\leq 2n}$ are both basis of 
$U_0^n\oplus U_0^n$. 
Let $A=(a_{ij})$ be the $2n\times 2n$ matrix given by
$$
\beta(Z_2,Z_j)=\sum_{r=1}^{2n}a_{rj}\beta(Z_1,Z_r).
$$
Let $\lambda\in\mathbb{C}$ be an eigenvalue of $A$ with
corresponding eigenvector $(v^1,\ldots,v^{2n})\in\mathbb{C}^{2n}$. 
Extending $\beta$ linearly from 
$V^{2n}\otimes\mathbb{C}$ to $W^{p,p}\otimes\mathbb{C}$, we have
$$
\sum_{j=1}^{2n}v^j\beta(Z_2,Z_j)
=\lambda\sum_{j=1}^{2n}v^j\beta(Z_1,Z_j).
$$
Hence $\beta(S,T)=0$ where $S=Z_2-\lambda Z_1$ and 
$T=\sum_{j=1}^{2n}v^jZ_j$. Setting $S=S_1+iS_2$ and 
$T=T_1+iT_2$, it follows that
$$
\beta(S_1,T_1)=\beta(S_2,T_2)\;\;\mbox{and}\;\;  
\beta(S_1,T_2)+\beta(S_2,T_1)=0.
$$
Then, if $X=S_1-JS_2$ and $Y=T_1+JT_2$ we obtain using 
\eqref{antihol} and \eqref{symmetries} that
\begin{align*}
\beta(X,Y)
&=\beta(S_1,T_1)+\beta(S_1,JT_2)-\beta(JS_2,T_1)-\beta(JS_2,JT_2)\\
&=\beta(S_1,T_1)+\beta(S_1,JT_2)+\beta(S_2,JT_1)-\beta(S_2,T_2)=0.
\end{align*}
Similarly, we obtain $\beta(X',Y')=0$ for $X'=S_1+JS_2$ and 
$Y'=T_1-JT_2$.  The vectors $X$ and $X'$ are both non-zero. 
For instance, if $X=0$ then 
$$
JS_2+iS_2=S=Z_2-\lambda Z_1.
$$
Thus $JS_2=Z_2-\text{Re}(\lambda)Z_1$ and $S_2=-\text{Im}(\lambda)Z_1$.
But then  $Z_2\in\spa\{Z_1,JZ_1\}$, and this is a contradiction.     
Finally, if we have that $Y=Y'=0$ then $T=0$, that is also a 
contradiction.  
\vspace{1ex}

\noindent{\it Fact} $2$.
There exists a vector $Z_0\in V^{2n}$ satisfying $\dim N(Z_0)=2n-2$. 
\vspace{1ex}

This Fact holds for $n=1$.  In fact, if $0\neq X\in V^2$ then
$\beta(X,X)$ and $\beta(X,JX)$ are linearly independent and hence
$N(X)=0$.  
\vspace{1ex}

We argue for $n\geq 2$. By Fact $1$ there are nonzero vectors
$X,Y\in V^{2n}$ such that $\beta(X,Y)=0$. From \eqref{symmetries} 
we have that $\dim B_XV=2r$ where $1\leq r\leq n-1$ since 
$B_XY=0$ and $X\notin N(X)$. Thus $\dim N(X)=2n-2r$.
Moreover, we may assume that $2\leq r\leq n-1$ since Fact $2$ 
holds when $r=1$ by taking $Z_0=X$.  Therefore, we assume that $n\geq 3$.

Set $U_1^t=\pi_1(B_X(V))$. We claim that $t=r$ and that
\be\label{betaximage}
B_X(V)=U_1^r\oplus U_1^r=\{\beta(Z,X)\colon Z\in V^{2n}\}.
\ee
In order to prove the claim, we first show that 
\be\label{UsUs}
B_X(V)+\{\beta(Z,X)\colon Z\in V^{2n}\}=U_1^t\oplus U_1^t
\ee
which, in particular, yields $t\geq r$. If 
$\beta(X,Z_1)=(\xi,\eta)$ and $\beta(Z_2,X)=(\zeta,\delta)$ 
then \eqref{symmetries} gives $\beta(X,JZ_1)=(\eta,-\xi)$ 
and $\beta(X,JZ_2)=-(\delta,\zeta)$. Hence 
$\xi,\zeta,\eta,\delta\in U_1^t$ and thus 
$(\xi,\eta),(\zeta,\delta)\in U_1^t\oplus U_1^t$. 
Now let $(\xi_1,\xi_2)\in U_1^t\oplus U_1^t$. Then there are 
$Z_1,Z_2\in V^{2p}$ such that $\beta(X,Z_i)=(\xi_i,\eta_i)$, 
$i=1,2$. Using \eqref{symmetries} we obtain 
$$
(\xi_1,\xi_2)=\frac{1}{2}\left(\beta(X,Z_1-JZ_2)+\beta(Z_1+JZ_2,X)\right),
$$
and \eqref{UsUs} has been proved. 

We show that $t=r$. In the process and for later use, 
we also prove that 
\be\label{cl}
\Sal(\bar\beta)=U_2^{n-r}\oplus U_2^{n-r}
\ee
where $\bar\beta=\beta|_{N(X)\times N(X)}$ and 
$U_2^{n-r}=U_0^n\cap(U_1^r)^\perp$. The flatness of $\beta$ 
gives
$$
\lp\beta(X,Z),\beta(S,T)\rp=\lp B_XT,\beta(S,Z)\rp=0
$$ 
for any $S,T\in N(X)$ and $Z\in V^{2n}$. From 
\eqref{symmetries} we have $\beta(S,X)=0$ and thus
$$
\lp\beta(Z,X),\beta(S,T)\rp=\lp\beta(Z,T),\beta(S,X)\rp=0
$$
for any $S,T\in N(X)$ and $Z\in V^{2n}$. Since 
$\Sal(\bar\beta)\subset\Sal(\beta)$ then from \eqref{UsUs} we have 
\be\label{tr}
\Sal(\bar\beta)\subset U_2^{n-t}\oplus U_2^{n-t}.
\ee

Suppose that the vector subspace $\Sal(\bar\beta)$ is 
nondegenerate. Since $\dim N(X)=2n-2r$ we have from
\eqref{estimnucleo2} and \eqref{tr} that
$$
\dim\mathcal{N}(\bar\beta)
\geq\dim N(X)-\dim\Sal(\bar\beta)\geq 2t-2r.
$$
Since $\mathcal{N}(\bar\beta)=0$ because 
$\beta(Z,Z)\neq 0$ if $Z\neq 0$ and  $t\geq r$, 
then  $t=r$ and $\dim\Sal(\bar\beta)=2n-2r$. 
Suppose now that the vector subspace $\Sal(\bar\beta)$ 
is degenerate.  If $n=3$ then $r=2$. Since 
$\Sal(\bar\beta)\neq 0$ then also $t=2$
and \eqref{cl} follows.  Then let $n\geq 4$. We have from
\eqref{tr} that
$$
\dim\Sal(\bar\beta)\leq 2n-2t\leq 2n-2r=\dim N(X) 
$$
and the assumption in the statement yields that $t=r$ and 
$\dim\Sal(\bar\beta)=2n-2r$. Having proved that 
$\dim\Sal(\bar\beta)=2n-2r$ holds in both 
cases then \eqref{cl} follows from \eqref{ebetaimage}. 
Finally, since $t=r$ then \eqref{betaximage} follows 
from \eqref{UsUs} using \eqref{symmetries}, and thus the 
claim has been proved.
\vspace{1ex}

We now claim that in fact the vector subspace $\Sal(\bar\beta)$ is 
nondegenerate. Assuming otherwise we have from 
\eqref{betadecomp} that 
\be\label{betarestriction}
\bar\beta(S,T)=\bar\beta_1(S,T)+2((S,T)v,(S,JT)v)
\ee 
for any $S,T\in N(X)$ where the light-like vector $v\in U_0^n$  
satisfies $\<v,w\>=-1$ and
\be\label{restrict}
\Sal(\bar\beta)\cap(\Sal(\bar\beta))^\perp
=\spa\{v\}\oplus\spa\{v\}.
\ee
Since $\beta(R,X)=0$ for  $R\in N(X)$ by \eqref{symmetries}, 
then the flatness of $\beta$ gives
$$
\lp\beta(Y,X),\beta(R,Z)\rp=\lp\beta(Y,Z),\beta(R,X)\rp=0
$$
for any $Y,Z\in V^{2n}$ and $R\in N(X)$. Therefore, for 
$R\in N(X)$ the vector subspaces $B_R(V)$ and 
$\{\beta(Y,X)\colon Y\in V^{2p}\}$ are orthogonal. 
From \eqref{betaximage} and \eqref{cl} we obtain that
$B_R(V)\subset U_2^{n-r}\oplus U_2^{n-r}
=\Sal(\bar\beta)$ if $R\in N(X)$.  
By \eqref{estimnucleo} the vector subspace 
$N_0=\mathcal{N}(\bar\beta_1)\subset N(X)$ 
is nonzero and from \eqref{symmetries} and \eqref{betarestriction} 
we have $B_S(N(X))=\spa\{v\}\oplus\spa\{v\}$ for any 
$S\in N_0$.  Then from \eqref{restrict} and the flatness of $\beta$ 
it follows that
$$
\lp\beta(S,Y),\beta(R,T)\rp=\lp B_ST,\beta(R,Y)\rp=0
$$
for any $S\in N_0$, $Y\in V^{2n}$ and $R,T\in N(X)$. Hence 
$B_S(V)\subset\Sal(\bar\beta)\cap(\Sal(\bar\beta))^\perp$ 
for any $S\in N_0$. Then from \eqref{betarestriction} and 
\eqref{restrict} it follows that 
$B_S(V)=\spa\{v\}\oplus\spa\{v\}$ for any $S\in N_0$. 
Since $\dim N(S)=2n-2$ then by Lemma \ref{kercod2} 
the subspace $B_S(V)$ should be nondegenerate. This
is a contradiction that proves the claim.

If $Z_0\in N(X)$ then \eqref{symmetries} gives that also
$\beta(Z_0,X)=0$. The flatness of $\beta$ yields
$$
\lp\beta(Z_0,Y),\beta(Z,X)\rp=\lp\beta(Z_0,X),\beta(Z,Y)\rp=0
$$
for any $Y,Z\in V^{2n}$. Then the second equality in 
\eqref{betaximage} gives $B_{Z_0}(V)\subset U_2^{n-r}\oplus U_2^{n-r}$. 
If $r=n-1$ then $Z_0$ satisfies $\dim\ker B_{Z_0}=2n-2$ as required. 
Hence, we assume that $r\leq n-2$ and thus $n\geq 4$ since $r\geq 2$.

We conclude the proof of Fact $2$ arguing by induction. 
Assume that it is true for any dimension until $n-1$. 
We have seen that $\dim N(X)=2n-2r$ with $2\leq r\leq n-2$. 
Thus the assumption of the induction applies to 
$\bar\beta\colon N(X)\times N(X)\to U_2^{n-r}\oplus U_2^{n-r}\subset W^{p,p}$ 
since the vector subspace  $\Sal(\bar\beta)=U_2^{n-r}\oplus U_2^{n-r}$ 
is nondegenerate. Hence there is $Z_0\in N(X)$ such that 
$\dim\ker\bar{B}_{Z_0}=2n-2r-2$. Then Lemma \ref{kercod2}
applies and by \eqref{kercod2decomp} and \eqref{cl} we have 
$$
U_2^{n-r}\oplus U_2^{n-r}
=\bar\beta_{Z_0}(N(X))\oplus\Sal(\bar\beta|_{\bar{N}(Z_0)\times\bar{N}(Z_0)})
$$
where $\bar{N}(Z_0)=\ker\bar{B}_{Z_0}$.  It follows that
$\dim\Sal(\bar\beta|_{\bar{N}(Z_0)\times\bar{N}(Z_0)})=2n-2r-2$.
By the flatness of $\beta$ we have 
$$
\lp\beta(Z_0,Y),\bar\beta(R,T)\rp=\lp\bar\beta(Z_0,T),\beta(R,Y)\rp=0
$$
for any $R,T\in\bar{N}(Z_0)$ and $Y\in V^{2n}$. Hence the vector
subspaces $\Sal(\bar\beta|_{\bar{N}(Z_0)\times \bar{N}(Z_0)})$ and
$B_{Z_0}(V)$ are orthogonal. Since 
$B_{Z_0}(V)\subset U_2^{n-r}\oplus U_2^{n-r}$ 
then $\dim B_{Z_0}(V)=2$, and hence Fact $2$ has been proved.
\vspace{1ex}

We conclude the proof of the proposition by means of a recursive 
construction. Notice that it suffices to construct an orthogonal 
basis of $\Sal(\beta)$ since it can be normalized. By Fact $2$ and 
\eqref{symmetries} there is $X_1\in V^{2n}$ such that the 
$J$-invariant vector subspace $N(X_1)$ satisfies $\dim N(X_1)=2n-2$. 
Moreover, $\beta(X_1,X_1)=(\xi_1,0)\neq 0$ and 
$\beta(X_1,JX_1)=(0,-\xi_1)$. By Lemma \ref{kercod2} the vector 
$\xi_1$ is either space-like or time-like.

By Lemma \ref{kercod2} the decomposition
$\Sal(\beta)=B_{X_1}(V)\oplus\Sal(\beta_{N(X_1)\times N(X_1)})$
is orthogonal and the vector subspace 
$\Sal(\beta_{N(X_1)\times N(X_1)})$ is nondegenerate and has 
dimension $2n-2$. 
By Fact $2$ there is $X_2\in N(X_1)$ such that 
$\dim\ker B_{X_2}\cap N(X_1)=2n-4$. Again we have that 
$\beta(X_2,X_2)=(\xi_2,0)\neq 0$ and 
$\beta(X_2,JX_2)=(0,-\xi_2)$ where $\xi_2$ is either 
space-like or time-like and the vectors $\xi_1$ and $\xi_2$ 
are orthogonal. Since $N(X_1)$ is $J$-invariant, then
$$
\beta(X_1,X_2)=0=\beta(X_1,JX_2)
$$
Then we have that $X_1,JX_1,X_2,JX_2$ are orthogonal. 
If we have $n=2$ then the desired basis is 
$X_1,JX_1,X_2,JX_2$ and if $n\geq 3$ we have to reiterate 
the construction.\qed

\section{The proofs}

This section is devoted to provide the proofs of the 
results stated in the introduction.
\vspace{1ex}

Let $f\colon M^{2n}\to\Hy^{2n+p}$ be an isometric immersion of a 
Kaehler manifold. Then let $g=i\circ f\colon M^{2n}\to\L^{2n+p+1}$
where $i\colon\Hy^{2n+p}\to\L^{2n+p+1}$ denotes the inclusion. We 
have that $N_gM=i_*N_fM\oplus\spa\{g\}$ and
$$
\a^g(X,Y)=i_*\a^f(X,Y)+\<X,Y\>g
$$
for any $X,Y\in\mathfrak{X}(M)$. At $x\in M^{2n}$ let
$\gamma,\beta\colon T_xM\times T_xM\to N_gM(x)\oplus N_gM(x)$
be the bilinear forms defined by
$$
\gamma(X,Y)=(\a^g(X,Y),\a^g(X,JY))\\
$$
and
\begin{align*}
\beta(X,Y)
&=\gamma(X,Y)+\gamma(JX,JY)\\
&=(\a^g(X,Y)+\a^g(JX,JY),\a^g(X,JY)-\a^g(JX,Y)).
\end{align*}
Then since $\a^g$ verifies 
the condition \eqref{shapeid} it follows that $\beta$ 
satisfies \eqref{nbeta}.

\begin{proposition}\po\label{flatforms}
Let $N_gM(x)\oplus N_gM(x)$ be endowed with the inner product
defined by
$$
\lp(\xi,\bar\xi),(\eta,\bar\eta)\rp
=\<\xi,\eta\>-\<\bar\xi,\bar\eta\>.
$$
Then the bilinear form $\beta$ is flat and we have that
\be\label{product*}
\lp\beta(X,Y),\gamma(Z,T)\rp=\lp\beta(X,T),\gamma(Z,Y)\rp
\ee
for any $X,Y,Z,T\in T_xM$.
\end{proposition}

\proof It is well-known that the curvature tensor of a Kaehler manifold 
$M^{2n}$ satisfies that $R(X,Y)JZ=JR(X,Y)Z$ for any $X,Y,Z\in T_xM$. 
Then this together with the Gauss equation for $g$ give
\begin{align*}
\lp&\gamma(X,T),\beta(Z,Y)\rp\\
&=\<\a(X,T),\a(Z,Y)\>+\<\a(X,T),\a(JZ,JY)\>-\<\a(X,JT),\a(Z,JY)\>\\
&\quad+\<\a(X,JT),\a(JZ,Y)\>\\
&=\<R(X,Z)Y,T\>+\<\a(X,Y),\a(Z,T)\>+\<R(X,JZ)JY,T\>\\
&\quad+\<\a(X,JY),\a(JZ,T)\>-\<R(X,Z)JY,JT\>-\<\a(X,JY),\a(Z,JT)\>\\
&\quad +\<R(X,JZ)Y,JT\>+\<\a(X,Y),\a(JZ,JT)\>\\
&=\<\a(X,Y),\a(Z,T)\>+\<\a(X,JY),\a(JZ,T)\>-\<\a(X,JY),\a(Z,JT)\>\\
&\quad+\<\a(X,Y),\a(JZ,JT)\>\\
&=\lp\gamma(X,Y),\beta(Z,T)\rp
\end{align*}
which proves \eqref{product*}.
Using \eqref{antihol} and \eqref{product*} we have
\begin{align*}
&\lp\gamma(JX,JY),\beta(Z,T)\rp=\lp\gamma(JX,T),\beta(Z,JY)\rp=
-\lp\gamma(JX,T),\beta(JZ,Y)\rp\\
&=-\lp\gamma(JX,Y),\beta(JZ,T)\rp=\lp\gamma(JX,Y),\beta(Z,JT)\rp
=\lp\gamma(JX,JT),\beta(Z,Y)\rp.
\end{align*}
Then
\begin{align*}
&\lp\beta(X,Y),\beta(Z,T)\rp=
\lp\gamma(X,Y),\beta(Z,T)\rp+\lp\gamma(JX,JY),\beta(Z,T)\rp\\
&=\lp\gamma(X,T),\beta(Z,Y)\rp+\lp\gamma(JX,JT),\beta(Z,Y)\rp
=\lp\beta(X,T),\beta(Z,Y)\rp,
\end{align*}
and this concludes the proof.\vspace{2ex}\qed

\noindent \emph{Proof of Theorem \ref{theorem}}: Observe that
$s\leq p+1\leq n-1$. At any point $M^{2n}$ the vector subspace 
$\Sal(\beta)$ is degenerate. In fact, if we have otherwise we obtain 
from \eqref{estimnucleo2} that $\dim\mathcal{N}(\beta)\geq 2n-2p-2\geq 2$, 
which is a contradiction. We have from Proposition \ref{betaimage} 
that $s\leq p$ and the proof now follows from part $(iii)$ of 
Proposition \ref{umbilical}.\qed 

\begin{remark}\label{r1}\po{\em For $f\colon M^{2n}\to\Hy^{2n+p}$ 
as in Theorem \ref{theorem} one can obtain from part $(i)$ of 
Proposition \ref{umbilical} that there is a smooth unit normal 
vector field $\eta\in\Gamma(N_fM)$ such that the second fundamental 
form satisfies $A_\eta^f=I$; see below the argument in the proof 
of Theorem \ref{m6}. 
An elementary proof gives that $f(M)$ is contained in a 
horosphere in $\Hy^{2n+p}$ if and only if $\eta$ is parallel 
in the normal connection.
}\end{remark}

\begin{lemma}\label{genclass}\po
Let $f\colon M^{2n}\to\Hy^{2n+p}$, $p\leq n-1$, be an isometric 
immersion of a Kaehler manifold. At any $y\in M^{2n}$ let
$\beta\colon T_yM\times T_yM\to W^{p+1,p+1}=N_gM(y)\oplus N_gM(y)$
be surjective. Moreover, assume that there is an orthogonal basis 
$\{X_j,JX_j\}_{1\leq j\leq n}$ of $T_yM$ such that parts 
$(i)$ and $(ii)$ in Proposition \ref{diag} hold. Then $p=n-1$ 
and $f(M)$ is an open subset  of the submanifold given by 
Example \ref{example}.
\end{lemma}

\proof  The vectors 
$(\xi_j,0)=\beta(X_j,X_j)\in N_gM(y)\oplus N_gM(y)$, $1\leq j\leq n$, 
in part $(ii)$ of Proposition \ref{diag} are an 
orthonormal basis of $W^{p+1,p+1}$. Thus $p=n-1$.
Part $(i)$ and \eqref{product*} give that
$$
\<A_{\xi_i}X_j,Y\>=\lp\gamma(Y,X_j),\beta(X_i,X_i)\rp
=\lp\gamma(Y,X_i),\beta(X_i,X_j)\rp=0\;\;\mbox{if}\;\;i\neq j,
$$
for any $Y\in T_yM$. Similarly, we have $\<A_{\xi_i}JX_j,Y\>=0$. 
Hence
\be\label{kernel}
\spa\{X_i,JX_i\}_{1\leq i\neq j\leq n}\subset\ker A_{\xi_j}
\ee
and thus the vector subspaces $\{\Ima A_{\xi_i}\}_{1\leq i\leq n}$ 
are orthogonal. Thus any pair of shape operator in  
$\{A_{\xi_j}\}_{1\leq j\leq n}$ commute, and we conclude from 
the Ricci equation that $g$ has flat normal bundle. 

From \eqref{kernel} we have $0\leq\rank A_{\xi_j}\leq 2$,
$1\leq j\leq n$. In fact, $\rank A_{\xi_j}=2$, 
$1\leq j\leq n$, for if $\rank A_{\xi_r}\leq 1$ for some 
$1\leq r\leq n$ then there is 
$0\neq Y\in\spa\{X_r,JX_r\}\cap\ker A_{\xi_r}$.  
But then \eqref{kernel} yields $Y\in\mathcal{N}(\a^g)$ in 
contradiction with $A_g=-I$.  

At any point $x\in M^{2n}$ there exist orthonormal normal vectors 
$\{\xi_j(x)\}_{1\leq j\leq n}$ and orthogonal tangent vectors 
$\{X_j(x), JX_j(x)\}_{1\leq j\leq n}$ such that 
$\ker A_{\xi_j(x)}=\oplus_{i\neq j}E_i(x)$  where 
we denote $E_j(x)=\spa\{X_j(x),JX_j(x)\}$.

From \eqref{kernel} and since $f$ has flat normal bundle, then
at any $x\in M^{2n}$ there is $r=r(x)$, $0\leq r\leq n$, such that 
the tangent space decomposes orthogonally as 
$$
T_xM=F_1(x)\oplus\cdots\oplus F_{n+r}(x),
$$
where $F_i(x)\oplus F_{n+i}(x)=E_i(x)$ if 
$0\leq i\leq r$ and $F_i(x)=E_i(x)$ if 
$r+1\leq i\leq n$, such that for each $\xi\in N_fM(x)$ there 
exist pairwise distinct $\lambda_i(\xi)\in\R$, $1\leq i\leq n+r$, 
satisfying
$$
A_\xi|_{F_i(x)}=\lambda_i(\xi)I.
$$
Since the maps $\xi\mapsto\lambda_i(\xi)$ are linear then there are 
unique pairwise distinct vectors $\eta_i(x)\in N_fM(x)$, $1\leq i\leq n+r$, 
such that
\be\label{shapediag}
\lambda_i(\xi)=\<\xi,\eta_i(x)\>,\; 1\leq i\leq n+r,
\ee
and
$$
F_i(x)=\{X\in T_xM\colon \a(X,Y)=\<X,Y\>\eta_i(x)
\;\;\mbox{for all}\;\ Y\in T_xM\}.
$$

On each connected component $M_r$ of the open dense subset of
$M^{2n}$ where the function $r(x)$ is constant, the maps 
$x\in M_r\mapsto\eta_i(x)$ and $x\in M_r\mapsto F_i(x)$, 
$1\leq i\leq n+r$, define smooth vectors fields and smooth 
vector bundles, respectively.   

If $F_i(x)\oplus F_{n+i}(x)=E_i(x)$ let $Y\in E_i(x)$ satisfy 
$A_{\xi_i}Y=\mu Y\neq 0$ and $A_{\xi_i}JY=\bar\mu JY\neq 0$. Let
$(Z)_F$ denote taking the
$F$-component of the vector $Z$.
Then \eqref{shapediag} gives
\begin{align*}
\mu Y
&=A_{\xi_i}Y
=A_{\xi_i}(Y)_{F_i}+A_{\xi_i}(Y)_{F_{n+i}}
=\<\xi_i,\eta_i\>(Y)_{F_i}+\<\xi_i,\eta_{n+i}\>(Y)_{F_{n+i}},\\
\bar\mu JY&=A_{\xi_i}JY
=A_{\xi_i}(JY)_{F_i}+A_{\xi_i}(JY)_{F_{n+i}}
=\<\xi_i,\eta_i\>(JY)_{F_i}+\<\xi_i,\eta_{n+i}\>(JY)_{F_{n+i}},\\
0&=A_{\xi_i}X_j=A_{\xi_i}(X_j)_{F_j}+A_{\xi_i}(X_j)_{F_{n+j}}
=\<\xi_i,\eta_j\>(X_j)_{F_j}+\<\xi_i,\eta_{n+j}\>(X_j)_{F_{n+j}},\\
0&=A_{\xi_i}JX_j=A_{\xi_i}(JX_j)_{F_j}+A_{\xi_i}(JX_j)_{F_{n+j}}
=\<\xi_i,\eta_j\>(JX_j)_{F_j}+\<\xi_i,\eta_{n+j}\>(JX_j)_{F_{n+j}}
\end{align*}
where $i\neq j$.
Thus the vectors $\eta_i$, $\eta_{n+i}$ are nonzero and satisfy 
$\eta_i,\eta_{n+i}\in\spa\{\xi_i\}$.  Hence the vectors fields
$\eta_1,\ldots,\eta_n$ form an orthogonal frame and 
$\xi_i\in\spa\{\eta_i\}$, $1\leq i\leq n$. It follows that on 
$M_r$ the orthonormal vector fields $\xi_1,\ldots,\xi_n$ are smooth.  

Next we argue on $M_r$. The Codazzi equation 
$(\nabla_SA)(T,\xi_j)=(\nabla_TA)(S,\xi_j)$ gives
$$
A_{\xi_j}[S,T]+A_{\nabla_S^\perp\xi_j}T=A_{\nabla_T^\perp\xi_j}S
$$
for any $S,T\in \ker A_{\xi_j}$.  Being the vector subspaces 
$\{\Ima A_{\xi_j}\}_{1\leq j \leq n}$ orthogonal, we have
$$
\<\nabla_S^\perp\xi_j,\xi_k\>A_{\xi_k}T
=\<\nabla_T^\perp\xi_j,\xi_k\>A_{\xi_k}S,\;\;
1\leq j\neq k\leq n
$$
for any $S,T\in\ker A_{\xi_j}$. Since $\rank A_{\xi_k}=2$  
then the $\xi_j$'s are parallel in the normal connection 
along $\ker A_{\xi_j}$, and hence along $TM_r$.

Let $Y\in E_j$ satisfy $A_{\xi_j}Y=\lambda Y$.
Then the  above Codazzi equation also gives
$$
\lambda\nabla_SY=A_{\xi_j}[S,Y]-S(\lambda)Y\in E_j
$$
for any $S\in\ker A_{\xi_j}$, $1\leq j\leq n$. Thus the  
subbundle $E_j$ is parallel along $\ker A_{\xi_j}$. Hence by 
\eqref{kernel} and since $\rank A_{\xi_j}=2$, we have 
 $$
\<\nabla_{X_j}Y,X_k\>=-\<Y,\nabla_{X_j}X_k\>=0,\;\;
\<\nabla_{JX_j}Y,X_k\>=-\<Y,\nabla_{JX_j}X_k\>=0,
$$
$$
\<\nabla_{X_j}Y,JX_k\>=-\<Y,\nabla_{X_j}JX_k\>=0, \;\;
\<\nabla_{JX_j}Y,JX_k\>=-\<Y,\nabla_{JX_j}JX_k\>=0
\;\;\mbox{if}\;\;j\neq k,
$$
and therefore the $E_j$'s are parallel along $M_r$. Then by the de Rham 
theorem there is an open subset $U\subset M_r$ which decomposes
in a Riemannian product $M_1^2\times\cdots\times M_n^2$ of surfaces. 
Since the codimension of $g$ is $n$ and $\a^g(Y_i,Y_j)=0$ if 
$Y_i\in E_i$ and $Y_j\in E_j$, $i\neq j$, 
then by Theorem $8.7$ in \cite{DT} there 
are isometric immersions $g_1\colon M_1^2\to\L^3$ and 
$g_j\colon M_j^2\to\R^3$, $2\leq j\leq n$, such that
$$
g|_U(x_1,\ldots,x_n)=(g_1(x_1),\ldots,g_n(x_n)).
$$
Since $g(M)\subset\Hy^{3n-1}$ then $\<g,g\>=-1$. Hence
$\<g_j{}_*X_j,g_j\>=\<g_*X_j,g\>=0$, and therefore 
$\|g_j\|=r_j$ with $-r_1^2+\sum_{j=2}^nr_j^2=-1$.
\vspace{2ex}\qed

We are now in conditions to give the remaining proofs of
the results stated in the introduction.
\vspace{2ex}

\noindent\emph{Proof of Theorem \ref{classification}}: 
First observe that it suffices to prove that this result holds  
on some open subset of $M^{2n}$. By Theorem \ref{theorem} we have
that $p=n-1$.  In an open neighborhood $U$ of $x_0$ in $M^{2n}$ 
there is a complex vector subbundle $\bar{V}\subset TM$ such that 
$\bar{V}(x_0)=V^{2m}$, $m\geq n-1$, and that $K(S,JS)>0$ for any 
$0\neq S\in\bar{V}$. Suppose that at $y\in U$ there is a $J$-invariant 
vector subspace $V_1\subset T_yU$ such that the subspace 
$\Sal(\beta|_{V_1\times V_1})$ is degenerate and 
$2\bar{s}=\dim\Sal(\beta|_{V_1\times V_1})\leq\dim V_1-2$. Then
by part $(iii)$ of Proposition \ref{umbilical} there exists a 
$J$-invariant vector $P^{2m}\subset T_yU$ with $m\geq n-\bar{s}+1$ 
such that $K(S,JS)\leq 0$ for any $S\in P^{2m}$.  Now by Proposition 
\ref{betaimage} we have $\bar{s}\leq p=n-1$. Hence $m\geq 2$ and this is a 
contradiction.  Thus the vector subspace $\Sal(\beta)$
is nondegenerate. From Proposition  \ref{diag} we obtain that 
$s=n=p+1$, that is, that $\beta$ is surjective. The proof now follows 
from Proposition~\ref{diag} and Lemma \ref{genclass}.
\qed\vspace{2ex}

\noindent\emph{Proof of Theorem \ref{unbounded}:} 
Recall that $\a^g(X,Y)=i_*\a^f(X,Y)+\<X,Y\>g$. 
Let $U$ be the subset of points of $M^{2n}$ such that at any $x\in U$ 
there is a $J$-invariant vector subspace $V_1\subset T_xM$ such that 
$\bar\beta=\beta|_{V_1\times V_1}$ satisfies that 
$\Sal(\bar\beta)$ is a degenerate subspace and  
$\dim\Sal(\bar\beta)\leq\dim V_1-2$.  If
the subspace $\Sal(\beta(x))$ degenerate for $x\in M^{2n}$ then
$x\in U$. In fact, we have for $s$ in \eqref{ebetaimage} that
$s<p+1$ since, otherwise,  $\Sal(\beta(x))=N_gM(x)\oplus N_gM(x)$ 
and therefore this subspace would be nondegenerate. Having 
that $p=n-1$, then $\dim\Sal(\beta)=2s\leq 2p=2n-2=\dim T_xM-2$, 
and thus $x\in U$.

Under the above conditions Proposition \ref{umbilical} applies 
to $\bar\beta(x)$ for $x\in U$. Then by part $(ii)$ there are 
a unit space-like vector $i_*\eta(x)\perp g(x)$ defined by 
$v=i_*\eta(x)+g(x)$ and the Lorentzian plane 
$L=\spa\{i_*\eta(x),g(x)\}$ such that for $P\subset N_fM(x)$ 
given by $i_*P=L^\perp$, we have 
\be\label{alphadecomumb}
\a^f(X,Y)(x)=\a^f_P(X,Y)(x)+\<X,Y\>\eta(x)
\ee
where $i_*\a^f_P(X,Y)=\a^g_{L^\perp}(X,Y)$. From \eqref{estimnucleo} 
applied to $\bar\beta$ we obtain that
$\bar\beta_1=\pi_{L^\perp\times L^\perp}\circ \bar\beta$ satisfies
$\dim\mathcal{N}(\bar\beta_1(x))\geq 4$. 
If $S\in\mathcal{N}(\bar\beta_1(x))$ then \eqref{betadecomp} gives
$$
0=\bar\beta_1(S,S)=(\a_{L^\perp}^g(S,S)+\a_{L^\perp}^g(JS,JS),0).
$$
It now follows from \eqref{alphadecomumb} that for any $x\in U$ we have
\be\label{alphadecomumb2}
\a^f(S,S)=\a^f_P(S,S)+\|S\|^2\eta(x)
\;\;\mbox{and}\;\;
\a^f(JS,JS)=-\a^f_P(S,S)+\|S\|^2\eta(x).
\ee

Suppose that there exists an open subset $U_0\subset M^{2n}$
such that $U\cap U_0=\emptyset$. We have by Proposition \ref{diag},
Lemma \ref{genclass} and Theorem \ref{classification} 
that $g=i\circ f\colon M^{2n}\to\L^{3n}$ is as 
given by the latter result. Then $g(M)$ is unbounded due to 
completeness of $M^{2n}$ and the presence of the hyperbolic 
factor, and hence also $f(M)$ is unbounded. 
\vspace{1ex}

By the above, it remains to consider the case when $U$ is 
dense in $M^{2n}$. We denote by  $r\colon\Hy^{3n-1}\to[0,+\infty)$ 
the distance function in $\Hy^{3n-1}$ to a reference point and 
set $h=\cosh(r\circ f)$. Then equation $(14)$ in \cite{DV} 
states that
\be\label{hessian}
\hess h(x)(X,X)=\cosh (r(f(x)))\|X\|^2
+\sinh(r(f(x)))\<\grad r(f(x)),\a^f(X,X)\>
\ee
for any $X\in T_xM$. On the other hand, we have from 
\eqref{alphadecomumb2} that
$$
\<\grad r(f(x)),\a^f(S,S)\>=\<\grad r(f(x)),\a^f_P(S,S)\>
+\<\grad r(f(x)),\eta(x)\>\|S\|^2
$$
and
$$
\<\grad r(f(x)),\a^f(JS,JS)\>=-\<\grad r(f(x)),\a^f_P(S,S)\>
+\<\grad r(f(x)),\eta(x)\>\|S\|^2
$$
for any $x\in U$ and $S\in\mathcal{N}(\bar\beta_1(x))$.  Then 
\eqref{hessian} yields
\begin{align}\label{hessianestim}
\hess h(x)(S,S)&+\hess h(x)(JS,JS)\\
&=2\cosh (r(f(x)))
+2\sinh (r(f(x)))\<\grad r(f(x)),\eta(x)\>\|S\|^2\nonumber
\end{align}
for any $x\in U$ and $S\in\mathcal{N}(\bar\beta_1(x))$. 

If $f(M)$ is bounded, by the Omori-Yau maximal principle 
for the Hessian there exists a sequence $\{x_k\}_{k\in\mathbb{N}}$ 
in $M^{2n}$ such that 
$$
\lim h(x_k)=\sup\{h\}<+\infty\;\;\mbox{and}\;\;
\hess h(x_k)(X_k,X_k)\leq \frac{1}{k}\|X_k\|^2\;\mbox{for any}\; 
X_k\in T_{x_k}M.
$$
Let $U_k\subset M^{2n}$ be an open neighborhood containing
$x_k$ such that 
$$
|h(x_k)-h(y)|<\frac{1}{k}
\;\;\mbox{and}\;\;
\hess h(y)(X,X)\leq \frac{2}{k}\|X\|^2
$$
for any $y\in U_k$ and $X\in T_yM$. By the above, we have 
that $U\cap U_k\neq\emptyset$ for any $k\in\mathbb{N}$.  
Then there exists a sequence 
$y_k\in U_k\cap U$ such that $r^*=\lim r(f(y_k))>0$ when 
$k\mapsto+\infty$. We obtain from  \eqref{hessianestim}
that 
$$
\frac{2}{k}\geq \cosh r(f(y_k)))
+\sinh r(f(y_k)))\<\grad r(f(y_k)),\eta(y_k)\>
$$
Taking the limit as $k\mapsto+\infty$ and using that 
$\|\grad r\|=1$ gives
$$
0\geq \cosh r^*-\sinh r^*>0,
$$
and this is a contradiction.  
\qed\vspace{2ex}

\noindent\emph{Proof of Theorem \ref{m6}:} 
Proposition \ref{diag} applies on an open subset $U\subset M^{2n}$ 
where $\Sal(\beta)$ is nondegenerate. We obtain that $n=3$ and that 
$g(U)$ is as described by Lemma \ref{genclass}. Then the 
result follows from Theorem \ref{classification}. Therefore, we may 
assume that $\Sal(\beta)$ is degenerate at any point of $M^{2n}$.
By \eqref{uimage} the subspaces $\Sal(\beta)\cap(\Sal(\beta))^\perp$ 
have dimension two and form a smooth vector subbundle.  
Hence, there is a smooth normal vector field $v\in\Gamma(N_gM)$ 
such that $\<v,g\>=-1$ and \eqref{umbilic} holds.  
Let $\eta\in\Gamma(N_fM)$ be the smooth unit vector field such that
$v=i_*\eta+g$. Then \eqref{umbilic} gives $A_\eta^f=I$.  
The Codazzi equation yields that $\eta$ is parallel in 
the normal connection along the open dense subset of nonflat
points of  $M^{2n}$, and then on all of $M^{2n}$.
Hence $f(M)$ is contained in a horosphere.\qed

\section*{Acknowledgment}

The second author thanks the Mathematics Department of the University
of Murcia where most of this work was developed for the kind 
hospitality during his visit.
\medskip

This research is part of the grant PID2021-124157NB-I00, funded by\\
MCIN/ AEI/10.13039/501100011033/ ``ERDF A way of making Europe"

\inputencoding{utf8}
{\renewcommand{\baselinestretch}{1}
\hspace*{-30ex}\begin{tabbing}
\indent \= Sergio J. Chion Aguirre \hspace{25ex} 
Marcos Dajczer\\
\>CENTRUM Catolica Graduate Business School, \hspace{0.5ex} 
IMPA -- Estrada Dona Castorina, 110\\
\> Pontificia Universidad Catolica del Perú, \hspace{8ex}
22460-320,\\
\>Lima, Perú,  \hspace{36.9ex} 
Rio de Janeiro -- Brazil\\
\> sjchiona@pucp.edu.pe,\hspace{27.1ex} 
marcos@impa.br
\end{tabbing}}

\begin{thebibliography}{lbl}
\bibitem{AMR} L. Al\'{i}as, P. Mastrolia and M. Rigoli, 
``Maximum Principles and Geometric Applications", 
Springer Monographs in Mathematics, 2016.

\bibitem{CD} M. do Carmo and M. Dajczer, 
\emph{Conformal Rigidity},  
Amer. J. Math. {\bf 109} (1987), 963--985.

\bibitem{CDa} S. Chion and M. Dajczer,
\emph{Minimal real Kaehler submanifolds},
Mat. Contemp. {\bf 49} (2022), 236--250.

\bibitem{DT} M. Dajczer and R. Tojeiro, 
``Submanifold theory. Beyond an introduction".
Series: Universitext. Springer, 2019.

\bibitem{DR} M. Dajczer and L. Rodr\'iguez, 
\emph{Rigidity of real Kaehler submanifolds}, 
Duke Math. J. {\bf 53} (1986),  211--220.

\bibitem{DV} M. Dajczer and Th. Vlachos,
\emph{Kaehler submanifolds of hyperbolic space},
Proc. Amer. Math. Soc. {\bf 148}, 4015--4024.

\bibitem{FHZ} L. Florit, W. Hui and F. Zheng,   
\emph{On real K\"ahler Euclidean submanifolds 
with non-negative Ricci curvature},
J. Eur. Math. Soc. {\bf 7} (2005), 1--11.

\bibitem{Ry}  P. Ryan,
\emph{K\"ahler manifolds as real hypersurfaces}, 
Duke Math. J. {\bf 40} (1973), 207--213.
\end{thebibliography}
\end{document}